\documentclass[10pt]{amsart}
\usepackage[T2A]{fontenc}
\usepackage[cp1251]{inputenc}
\usepackage[english]{babel}
\usepackage{amsmath,latexsym,amsthm,amsfonts,tipa,upgreek}
\usepackage{amssymb,amscd,graphpap, stmaryrd}

\begin{document}

\begin{center}
\Large{Locally uniformly rotund renormings  of the spaces 
of continuous functions on  Fedorchuk compacts}
\end{center}

\bigskip
\begin{center}
S.P.Gul'ko${}^a$, A.V.Ivanov${}^b$, M.S.Shulikina${}^c$, S.Troyanski${}^d$\footnote{First three authors was 
 supported by the Russian Foundation for Basic Research in the framework of 
the scientific project N 17-51-18051. The second author was 
also supported by state order to the Karelian Research Centre 
of the Russian Academy of Sciences (Institute of Applied 
Mathematical Research KarRC RAS). The forth author was 
partially supported by MTM 2017-86182-P (AEI/FEDER,UE)
 and by Bulgarian National Scientific Fund under Grant 
DFNI/Russia, 01/06/23.06.2017.}
\end{center}

\bigskip
${}^a$ Faculty of Mechanics and Mathematics, Tomsk State 
University, Prospect Lenina 36 634010 Tomsk, Russian 
Federation E-mail address: gulko@math.tsu.ru

${}^b$ Institute of Applied Mathematical Research of 
Karelian Research Centre, Russian Academy of 
Sciences, Russian Federation,\\ E-mail address: 
alvlivanov@krc.karelia.ru 

${}^c$ Department of Mathematics and Computer Science, Tomsk Polytechnic University, Russian Federation, E-mail address:shulikinams@tpu.ru

${}^d$ Institute of Mathematics and Informatics, Bulgarian 
Academy of Science, bl.8, acad. G. Bonchev str., 1113 
Sofia, Bulgaria and 
Departamento de Matem\'aticas, Universidad de Murcia, 
Campus de Espinardo, 30100 Murcia, Spain, E-mail 
address:stroya@um.es

\bigskip

\begin{center}
{\bf Abstract}
\end{center}
\bigskip

We show that $C(X)$ admits an equivalent pointwise lower 
semicontinuous locally uniformly rotund norm provided
 $X$ is 
Fedorchuk compact of spectral height 3. In other words  $X$
 admits a fully closed map $f$
 onto a metric compact $Y$  such that $f^{-1}(y)$
 is metrizable for all  $y\in Y$.
  A continuous map of compacts $f:X\to Y$
is said to be fully closed if for any disjoint closed subsets
$A,B\subset X$  the intersection $f(A)\cap f(B)$
is finite. For instance the projection of the lexicographic 
square  onto the first factor is fully closed and all its fibers are 
homeomorphic to the closed interval.
\bigskip

\section{INTRODUCTION}

Let us recall that a Banach space $E$
(or the norm in $E$) is said to be locally uniformly rotund 
(LUR for short) if $\lim_n||x_n-x||=0$
  whenever 
  $$
\lim_n||(x_n+x)/2||=\lim_n||x_n||=||x||.
$$
The spaces with this property are at the core of renorming 
theory in Banach spaces and consequently have been 
extensively studied (see, for example, [3] and 
[14] and its reference). It is well known that the 
spaces with a LUR norm have the Kadec property, that is 
on the unit sphere of $E$ the norm topology and the weak 
topology coincide.  The LUR renorming techniques for a 
Banach space developed until now, are based in two 
different approaches. In the first one, for enough convex 
functions on the Banach space are constructed to apply 
Deville’s master lemma (see the decomposition method 
[3, Chapter 7, Lemma 1.1p.279]), sometimes adding 
an iteration processes and Banach’s contraction mapping 
theorem, to finally get an equivalent LUR norm. 
Originally this method use  the powerful method of 
projectional resolutions of the identity(PRI for short) The 
second one  is based to
a characterization of those Banach 
spaces that admit a LUR renorming by means of a linear 
topological condition of covering 
type [13]. The existence of 
such norm is deduced from the existence of some  maps 
acting from  normed space $E$ 
into metric space $Y$. This maps admit some covering properties (see [14]).  We present a particular result 
 in this 
direction which we use latter. This is nonlinear 
(convex) version for LUR renorming of Banach 
spaces with strong Markushevich basis.

{\bf Theorem 1.1.} {\it Assume that there is a locally bounded map 
$\Phi$   from a normed space $E$ into $c_0(\Gamma)$
for some set $\Gamma$ such that:

(i) for every finite set $A\subset \Gamma$
 is specified a separable subspace $Z_A$
 of $E$  such that:
 
(a) $Z_A\subset Z_B$ whenever $A\subset B\subset\Gamma$;

(b) $x\in \overline{span\cup_{n\in N}Z_{K_n}}^{||\cdot||}$  whenever  $x\in E$ and 
$\{K_n:n\in N\}$ is an increasing sequence (i.e. $K_n\subset K_m$ for $n<m$) of finite subsets of $\Gamma$ with $\cup_{n\in N}K_n\supseteq
supp(\Phi x)=\{\gamma\in\Gamma:\Phi x(\gamma)\not=0\}$;

(ii) there exists norming subspace  $F$ of $E^*$ such that for 
every $\gamma\in\Gamma$  the real function $\delta_\gamma\circ\Phi$
 on $E$  is non-negative, convex and $\sigma(E,F)$-lower semi-continuous, where $\delta_\gamma$
is the Dirac measure on $\Gamma$
at $\gamma$.

                                                                                                                                                               Then $E$ admits an equivalent  $\sigma(E,F)$-lower semi-continuous
 LUR norm.}
 
{\bf Remark.}
Condition(b) is equivalent to  

(b') for every $x\in E$ there exists an increasing sequence
$\{K_n:n\in N\}$  of finite subsets of
the support  $supp(\Phi x)$  with $x\in \overline{span\cup_{n\in N}Z_{K_n}}^{||\cdot||}$.   
 
{\bf Proof.} Let us consider the condition

(iii) for every $x\in E$ there exists separable subspace 
$Z(x)$ of $E$ with 
$$
x\in\overline{span\cup\{Z(x_n):n\in N\}}^{||\cdot||} \eqno(1.1)
$$
whenever $\{x_n:n\in N\}$
is a  sequence in $E$
 with $\lim_n||\Phi x_n-\Phi x||_\infty=0$.
 
In [14, Theorem 1.15, Theorem 2.14, Corollary 
1.21, Theorem 3.28, Corollary 4.34] is 
shown that
(ii) and (iii) imply existence in $E$
an equivalent $\sigma(E,F)$-lower semi-continuous LUR norm. 
So in order to prove the former Theorem we have to 
show (i) $=>$ (iii).

Let us mention that  from (a) it follows that
$$
span \cup_{n\in N}Z_{A_n}=span \cup_{n\in N}Z_{B_n}
$$
for any increasing sequences  $\{A_n:n\in N\}$, $\{B_n:n\in N\}$ 
of finite subsets of $\Gamma$ with $\cup_{n\in N}A_n=\cup_{n\in N}B_n$.
So we can define $Z_A$ 
for any countable subset $A$  of $\Gamma$ by
setting 
$$
Z_A=\overline{span\cup_{n\in N}Z_{A_n}}^{||\cdot||} 
$$
if $A=\cup_{n\in N} A_n$ where $\{A_n:n\in N\}$ an increasing sequence
of finite subsets of $\Gamma$. Clearly (a) holds for countable sets as well.

Set $Z(x)=Z_{supp(\Phi x)}$. Let $\{x_n:n\in N\}$
be a sequence in $E$ with                
$$
\lim_n||\Phi x_n-\Phi x||_\infty=0 \eqno(1.2)
$$
for some $x\in E$. Using (b') we can find an increasing sequence $\{K_m:m\in N\}$   of 
finite subsets of $supp(\Phi x)$
such that
$$
x\in \overline{span\cup_{m\in N}Z_{K_m}}^{||\cdot||} \eqno(1.3)
$$

Set $B_n=supp(\Phi x_n)$. Fix $m\in N$.
From (1.2) it follows that $K_m\subset B_n$
  for all large enough  $n$.
 Hence $Z_{K_m}\subset Z_{B_n}=Z(x_n)$ for all 
large enough  $n$. So
$$
\cup_{m\in N}Z_{K_m}\subset\cup_{n\in N}Z(x_n).
$$
This and (1.3) imply(1.1). $\Box$ 

{\bf Example 1.2.} Assume that $E$ 
has a srtrong Markushevich basis $\{e_\gamma:\gamma\in \Gamma\}$
 with conjugate system  $\{e^*_\gamma:\gamma\in\Gamma\}$, that is $e^*_\beta(e_\gamma)=\delta_{\beta\gamma}$, and for 
every $x\in E$
we have
$$
x\in\overline{span\{e_\gamma:e^*_\gamma(x)\not=0,\gamma\in\Gamma\}}^{||\cdot||}.
$$
For finite set $A\subset \Gamma$ set $Z_A=span\{e_\gamma:\gamma\in A\}$ and
define $\Phi:E\to c_0(\Gamma)$
by formula   $\Phi x(\gamma)=|e^*_\gamma(x)|.$       
From the definition of strong 
Markushevich basis for every $x\in E$
we have 
$$
x\in\overline{span\{e_\gamma:\Phi x(\gamma)=|e^*_\gamma(x)|\not=0,\gamma\in\Gamma\}}^{||\cdot||}.
$$

So $\Phi$ satisfies the hypothesis of  the former 
Theorem.

When Banach space $E$ admits PRI in the similar way 
we can construct map  $\Phi$. 
 
Lot of papers are devoted to find different classes of Hausdorff
compact spaces (compacts) $X $
for which $C(X)$ admits an equivalent 
pointwise lower semi-continuous LUR norm. Moving on 
now to topological properties, we say that a compact space 
$X$ is {\it Eberlein}  
if it is homeomorphic to a weakly compact subset of a Banach space. Equivalently, $X$
 is Eberlein if and only if $C(X)$ is a 
{\it weakly compactly generated} 
space. The space $X$  is called {\it Talagrand} or 
{\it Gul’ko }
compact if $C(X)$ is 
{\it weakly $K$-analytic} or {\it weakly countably 
determined},  respectively.
 For a full treatment of these concepts, we refer the reader to e.g.,  [3]. It turns out that if $X$ 
hails from one of these three classes of compact 
spaces then it can be treated as a subset of a pointwise 
compact cube $[0, 1]^\Gamma$, in such a way that, given
 $t \in X$, its coordinates $t(\gamma),\ \gamma\in\Gamma$
 behave according to certain rules. For instance, 
$X$ is Eberlein if and only if we can find such 
an embedding with the property that for any $t\in X$ 
 and $\varepsilon>0$, there are only finitely many 
$\gamma\in\Gamma$  satisfying $|t(\gamma)|\geq\varepsilon$.  
Similarly, spaces from the three corresponding 
classes of Banach spaces may be endowed with 
‘coordinate systems’, which permit the spaces to be 
carefully  analyzed in ways which are not feasible in a 
fully general, non-separable setting. Every compact space  
$X$ from the classes above shares the property that $X$
 may be 
embedded in $[0, 1]^\Gamma$, such that given any $t \in X$, the 
 support
 of $t$  is countable.  
 In general, a space satisfying this property is called 
{\it Corson} compact. 
The space $X$ is called {\it Valdivia} compact if it is as above, 
but, in this case, only a pointwise dense subset of points of 
$X\subset[0,1]^\Gamma$ 
are required to have countable support.  These 
classes have long been relevant to renorming theory. For 
example, it can be shown that if $X$  is Valdivia compact 
then $C(X)$admits a pointwise lower semi-continuous LUR 
norm (cf. [3, Corollary VII.1.10]) . In all this cases 
$C(X)$ admits PRI. Of course all this result can be obtained 
directly from Theorem 1.1.
 We mention that in $C(X)$
 spaces we have a canonical map to $c_0(\Gamma)$ 
 (see e.g. [14]). Indeed if $X\subset[0,1]^\Gamma$  
the uniform continuity of every $h\in C(X)$
 allows us to define the oscillation map $\Omega:C(X)\to c_0(\Gamma)$
  by formula 
$$
\Omega(h)(\gamma)=\sup\{h(t)-h(s):t,s\in X,(t-s)|_{\Gamma\setminus\{\gamma\}}=0\},
$$
where $h\in C(X)$.
Map $\Omega$ was introduced(see [4]) looking 
for countable sets of coordinates which control a 
continuous function to obtain extensions of the theorem of 
Mibu. 
It is easy to see that $\Omega$
is a bounded map  satisfying 
condition (ii) of Theorem 1.1. In [14, 2.7] is shown that 
$\Omega$  satisfies condition (iii) of Theorem 1.1 when 
$X$  is Helly compact of monotone functions on [0,1]. In this way 
is proved that $C(X)$
is LUR renormable when $X$ 
is Helly compact. A generalization of this result when 
$X$ is a particular case of Rosenthal compacts  can be 
found in[10], see also [12].

The aim of this note is to find a new class of compact 
spaces for which the corresponding space of continuous 
functions is LUR renormable.

The following definition goes back to V.V.Fedorchuk (see e.g.
[5, II.1.6] )

{\bf Definition 1.3.} {\it Continuous  map of compacts $f:X\to Y$
 is said to be fully 
closed if for every disjoint closed subsets $A$
 and $B$ of $X$ the set $f(A)\cap f(B)$
 is finite.}
 
Our  main result is next

{\bf Theorem 1.4.} {\it Let $X$
 be a compact space admitting a fully 
closed map $f$ onto a metrizable compact $Y$
such that the fibers $f^{-1}(y)$
are  metrizable  for every $y\in Y$.
 Then $C(X)$
admits an equivalent 
pointwise lower semi-continuous
LUR norm.}

The above class of compact spaces  is particular case of 
Fedorchuk compacts. The section 2 is devoted to this 
class. Now we give an example. Denote with
 $L$  the 
lexicographic square (the projection of this square onto the 
first factor is fully closed and all its fibers are 
homeomorphic to the closed interval). 
In [1] is shown that $C(L)$ is LUR renormable(for general 
case of totally ordered compacts see [9] and 
[14]). 

Let $f:X\to Y$  be a continuous map of compacts. Given $y\in Y$   define   
$$
osc_{f^{-1}(y)}(h)=\sup_{s,t\in f^{-1}(y)}\{h(t)-h(s)\}=diam(h(f^{-1}(y))
$$
for $h\in C(X)$. Clearly $osc_{f^{-1}(y)}(\cdot)$
is a pointwise lower semi-continuous semi-norm in 
$C(X)$ and $osc_{f^{-1}(y)}(\cdot)\leq2||\cdot||_\infty$.

We introduce fiberwise oscillation map $\Omega_f:C(X)\to l_\infty(Y)$
by formula 
$$
\Omega_f(h)(y)=osc_{f^{-1}(y)}(h)
$$
where $h\in C(X)$. 

Since $||\Omega_f(h)||_\infty\leq 2||h||_\infty$
 for all $h\in C(X)$
 the map $\Omega_f$  is bounded.  Since 
$osc_{f^{-1}(y)}(\cdot)$ 
is a pointwise lower semi-continuous semi-norm in 
$C(X)$ and                                                                                                                                                                                                                                                                                                                            
$\delta_y\circ\Omega_f(h)=osc_{f^{-1}(y)}(h)$ we get that  
$\Omega_f$  
satisfies condition (ii) of Theorem 1.1.
  In sections 2 and 3 we show that if $f$ a fully closed map satisfying the conditions of Theorem 1.4 then $\Omega_f$
maps $C(X)$
 into $c_0(Y)$
 and satisfies condition (i) of Theorem 1.1.

\section{FEDORCHUK COMPACTS}

The class of  Fedorchuk compact spaces was defined in 1984 [11] with the purpose of clarifying the limits of the application of the method of resolutions, which was developed by  Fedorchuk (see [5]) and showed exceptional efficiency in constructing counterexamples in general topology.
The definition of this class is inextricably linked with the concept of a fully closed mapping introduced by Fedorchuk in the process of developing the method mentioned above.  The original definition of fully closed mapping  was  cumbersome.
Later, his author obtained a number of  unobvious equivalent formulations (see [5, II. 1.6]).  The shortest of them for mappings of Hausdorff compact  spaces is the Definition 1.3 above.

{\bf Definition 2.1.} {\it Hausdorff compact  space  $ X $ is called a Fedorchuk compact 
  (or an $F$-compact) if there exists a well-ordered continuous inverse system (an $F$-system)
 $ S = \{X_\alpha, \pi^\alpha_\beta: \alpha, \beta \in \gamma \} $ (here $ \gamma $ is an ordinal number) giving in the limit $X$, in which $X_0$  is the point, all the neighboring projections $ \pi^{\alpha+1}_\alpha $ are fully closed ($ \alpha + 1 < \gamma $), and the inverse images of the points $ (\pi^{\alpha + 1}_\alpha)^{- 1}(x) $ are metrizable for any $ x \in X_\alpha $.  The spectral height $ sh(X) $ of an $F$-compact $ X $ is the smallest possible length $ \gamma $ of such a system.}
 
We consider only the Fedorchuk compacts  of spectral height 3. If $ X $ is an $F$-compact and $ sh(X) = 3 $, then $X$ is the limit of an $F$-system consisting of three spaces: $X_0, \  X_1 $ and $ X_2 $.  The limit of this system coincides with $ X_2 $, from which it follows that the compact space $ X = X_2 $ is non-metrizable, since otherwise it could be obtained as the limit of the $F$-system of two compacts: $ X_0 $ and $ X_2 $ (the map to a point is always fully closed).
Consequently, for any Fedorchuk compact of spectral height 3 there exists a fully closed map $ \pi^2_1: X = X_2 \to X_1 $
 onto a metric compact space for which the sets $ (\pi^2_1)^{- 1} (x), \ x \in X_1 $ are metrizable.
 
 Thus, $F$-compacts of spectral height 3 can be characterized as non-metrizable compacts that admit a fully closed map onto a metric compact with metrizable fibers.  Note that if such a fully closed map for a compact $ X $ exists, then it is almost unique in the following sense: almost all (that is, all but perhaps a countable set) nontrivial fibers of any two such maps coincide (see [6]). Remark, that 
only a point is an $F$-compact of  spectral height 1 and $F$-compacts of spectral height 2 are all non-trivial metric compact spaces. 

As noted in the introduction, the lexicographic square $L$ is an $F$-compact of spectral height 3.
The classical  space "two arrows" (the lexicographic product of unit segment and two-point set $[0,1]\times\{0,1\}$) also is $F$-compact of  spectral height 3 (the standard projection of this space onto the segment is a fully closed map) as well as the space "Alexandroff double circle"   (the projection of "double circle" onto the circle is fully closed).

To the class of $ F$-compacts of spectral height 3 belongs a series of various counterexamples constructed by  Fedorchuk.  Here, first of all, we should mention a group of compacts with non-coinciding dimensions.  Among them is the famous example of a two-dimensional compact that  has no partitions of lower dimension (see [2], page 314), and also a homogeneous separable compact $ X $ with the first axiom of countability with $ dim (X) = 1 <ind (X ) = 2 $ and a perfectly normal compact (constructed under the assumption of $CH$) with same values of  the dimensions $ dim $ and $ ind $ (see [5, III. 3.6 and 3.10]).

We will need the following characteristic property of fully closed map obtained by Fedorchuk [5]. 
Let $ f: X \to Y $ be a continuous map, and $A$ be an arbitrary subset of $Y$.  Consider a partition of a compact $X$ whose nontrivial elements are  sets $ f^{-1} (y) $ for $ y \in Y \setminus A $.  Let $ Y_A $ be the quotient space corresponding to this partition (with respect to 
$f$), that is
$$
Y_A=\{f^{-1}(y):y\in Y \setminus A\}\cup\{\{x\}:f(x)\in A\}.
$$
Let
 $ f_A: X \to Y_A$ be a    quotient  map, and let $\pi_A: Y_A\to Y$ be the unique map for which $f=\pi_A \circ f_A$. 

{\bf Proposition 2.2.} [5, II.  1.6].  {\it The map $ f $ is fully closed if and only if for each $ A \subset Y $ the space $ Y_A $ is Hausdorff.}

In [5, II, 1.7, 1.10] it is also shown  that for a fully closed map $ f $ the maps $f_A$ and $\pi_A $ are also fully closed for any $ A \subset Y$.

{\bf Proposition 2.3.} [5, II.  3.10].  {\it Let $ f: X \to Y $ be a fully closed map of the compact $ X $ onto a metric compact $ Y $ with metrizable fibers $ f^{-1}(y)$, $ y \in Y $.  The compact $X$ is metrizable if and only if the set of nontrivial fibers (i.e. the fibers 
$f^{-1}(y)$ that is not a 
singleton) of $f$ is countable.}

\section{PROOF OF THEOREM 1.4.}

The following assertion gives a characterization of fully closed maps, which plays an important role in what follows.

 {\bf Proposition 3.1.} {\it The map $ f: X \to Y $ of compacts is fully closed if and only if for any continuous map $ g: X \to K $ into a metric compact space $ K $ and any $ \varepsilon> 0 $ the set $ H_{g, \varepsilon} = \{y\in Y: diam(g(f^{- 1}(y))) \geq \varepsilon \} $ is finite.}
 
{\bf Proof.} Necessity.  Let $f$ be fully closed.  This means that for any disjoint closed subsets $ A, B \subset X $ the intersection $ f(A)\cap f(B)$ is finite.  Suppose that there exists a continuous map $ g: X \to K $ to the metric compact $K$ and $ \varepsilon> 0 $ such that $H_{g,\varepsilon}$ is infinite.  Consider the infinite family of closed subsets of $K$ indexed by   $y$:
$$
G=\{g(f^{-1}(y)):y\in H_{g,\varepsilon}\}.
$$
Let $ F $ be an    accumulation point of the family $ G $ in the space $ \exp (K) $ of nonempty closed subsets of $K$ endowed with the Hausdorff metric.  It is obvious that $ diam(F) \geq \varepsilon $.  We take in $ F $ two distinct points $ x_1, x_2 $ and their neighborhoods $ O_{x_1}, O_{x_2} $ with disjoint closures.  By the choice of $ F $, there exists an infinite subset $ D \subset H_{g, \varepsilon}$ such that $ g(f^{- 1}(y)) \cap O_{x_i} \not = \emptyset $ for any $ y  \in D $ for $ i = 1,2 $.  Consider disjoint closed subsets $ A_i= g^{- 1}\overline{O_{x_i}} $, $ i = 1,2 $ in the space $ X $.  We have $ D \subset f(A_1) \cap f(A_2) $, which contradicts that $f$ is fully closed.

Sufficiency. 
 Let $ f $ be not fully closed, that is, there exist  disjoint closed subsets $ A, B $ in $ X $ such that the intersection $ E = f (A) \cap f(B)$ is infinite.  Let $ g: X \to [0,1] $ be a continuous function on $ X $ that separates $A$ and $B$.  Then the set $ H_{g, 1} $ contains $ E $ and, consequently, is infinite.  $\Box$
 
{\bf Corollary 3.2.} {\it For any fully closed map of compacts $f:X\to Y$
 the fiberwise oscillation map $\Omega_f$ maps $C(X)$
 into $c_0(Y)$.}
 
Let $ f: X \to Y $ be a fully closed map of compacts. Given $A\subset Y$ define
$$
Z_A=\{h\in C(X):supp(\Omega_f(h))\subset A\}. \eqno(3.1)
$$
From the definition of $\Omega_f$ it
follows that $Z_A$ is  a 
closed subspace of  $C(X)$.

{\bf Lemma 3.3.} {\it  The spaces $Z_A$
 and $C(Y_A)$
  are isomorphically isometric. The linear operator
$T_A:Z_A\to C(Y_A)$
 defined by formula  
 $$
T_Ah=h\circ f_A^{-1} \eqno(3.2)
$$
give the isometry , i.e.
$$
||T_Ah||_\infty=||h||_\infty. \eqno(3.3)
$$                           
Moreover for every $y\in A$ and $h\in Z_A$
we have
$$
T_Ah(\pi_A^{-1}(y))=h(f^{-1}(y)). \eqno(3.4)
$$}

{\bf Proof.} We have
$$
Y_A=\{f^{-1}(y):y\in Y\setminus A\}\cup \{\{x\}:f(x)\in A\}.
$$
Pick $h\in Z_A$. Then $osc_{f^{-1}(y)}(h)=0$ if $y\in Y\setminus A$.
So $(T_Ah)(f^{-1}(y))=h\circ f_A^{-1}(f^{-1}(y))=h(x)$ (here $f^{-1}(y)$ is a point of $Y_A$) for every $x\in f^{-1}(y)$ and $y\in Y\setminus A$. If $f(x)\in A$ then $(T_Ah)(\{x\})=h(f_A^{-1}(\{x\}))=h(x)$.

This implies that $T_Ah$
 is a real single valued and continuous function on $Y_A$. 
Moreover (3.3) holds.

Since $f=\pi_A\circ f_A$ we get $f^{-1}(y)=f_A^{-1}(\pi_A^{-1}(y))$. 
So for $h\in Z_A$ and $y\in A$
$$
h(f^{-1}(y))=h(f_A^{-1}(\pi_A^{-1}(y)))=T_Ah(\pi_A^{-1}(y)).
$$
$\Box$

The main result of this section is

{\bf Proposition 3.4.} {\it Let $ f: X \to Y $ be a fully closed map of compact $X$ onto  a metric compact $Y$ with metrizable fibers $f^{-1} (y), y \in Y $. Let $K$ be a finite subset of a
countable set $A\subset Y$.   Then for every  $h\in Z_A$
 we have
 $$ 
dist(h,Z_K)=\inf \{||h-g||_\infty:g\in Z_K\}\leq||\Omega_f(h)|_{A\setminus K}||_\infty. \eqno(3.5)
$$ 
  }
  
The proof is based on the properties of 
Fedorchuk compacts. For this reason at first, before the 
proof, we introduce some concepts and notations. In the subsequent arguments we need the concept of a small image of a set.  Recall that  a small image $f^\sharp(B)$ of the set $B\subset X$  under a map $f:X\to Y$ is  the set $f^\sharp(B)=Y\setminus f(X\setminus B)$. Clearly $f^\sharp(B)\subset f(B)$. 
Let us mention  that for any continuous map $f:X\to Y$ of compacts and any open subset $U\subset X$, a small image $f^\sharp(U)$ is open in Y, since the set $f(X\setminus U)$ is closed in $Y$ as a continuous image of a compact.

We denote with $Fr(B)$ the topological  boundary of 
the set $B$, that is  $ Fr(B) =\overline{B}\setminus int(B)$.

Using the notations of section 2 we consider the 
fully closed map $\pi_A:Y_A\to Y$.
 Consider a partition of compact $Y_A$
 whose nontrivial elements are the sets $\pi_A^{-1}(y)$ for $y\in Y\setminus K$. Let $Y_{A,K}$
be the quotient space corresponding to this partition,  let   
$\pi_{A,K}:Y_A\to Y_{A,K}$
 be  the quotient map, and let  $\varkappa_{A,K}:Y_{A,K}\to Y$
 be the unique map for which $\pi_A=\varkappa_{A,K}\circ\pi_{A,K}$.
By Proposition 2.3 we have that $Y_{A,K}$
is metrizable since the set of nontrivial fibers $\{\varkappa_{A,K}^{-1}(y):\ y\in K\}$ of the fully closed map $\varkappa_{A,K}$ is finite and all this fibers and $Y$ are metrizable. We fix a metric $d$ on $Y_{A,K}$
that is compatible with the topology.  We denote 
with $B(u,\rho)=\{y\in Y_{A,K}:d(u,y)<\rho\}$  
the open ball centered at $u$ with  radius $\rho>0$. 

{\bf Proof of Proposition 3.4.} Pick $h\in Z_A$
 and set          
$$
s=||\Omega_f(h)|_{A\setminus K}||_\infty,\ g=T_Ah.
$$
So $g=h\circ f^{-1}_A\in C(Y_A)$. Set $D=A\setminus K$.
Enumerate the points $D$ by the natural numbers:
 $D=\{y_n:n\in N\}$ 
and denote by $u_n$
 the unique point in $Y_{A,K}$
with $\varkappa_{A,K}(u_n)=y_n$.
Set  $E=\{u_n:n\in N\}$.             

Fix $c>1$. For every $n\in N$
choose an open neighborhood $O_{y_n}$ of $\pi_A^{-1}(y_n)$ in $Y_A$ 
in such a way that
$$
diam(g(O_{y_n}))<c\cdot diam(g(\pi_A^{-1}(y_n)). \eqno(3.6)
$$
From (3.4) we get $diam(g(\pi_A^{-1}(y_n)))=diam(h(f^{-1}(y_n)))$.
Since  $y_n\in A\setminus K$
we get  $diam(h(f^{-1}(y_n)))\leq s$.
 So 
$$
diam(g(O_{y_n}))<cs. \eqno(3.7)
$$

{\bf Claim 1.} {\it There exists a sequence of 
open sets  $\{U_n:n\in N\}$  in $Y_{A,K}$
such that $E\subset\cup_{n\in N}U_n$,
the closures of  $U_n$ are pairwise disjoint,  $\overline{U_n}\subset\pi^\sharp_{A,K}(O_{y_n}),\ Fr(U_n)\cap E=\emptyset$ and $diam(\overline{U_n})<\min\{cs,1/n\}$.}   

{\bf Proof of Claim 1.} We construct the open sets $U_n,\ n\in N$
by recursion as follows.

Step 1. Clearly $u_1\in \pi^\sharp_{A,K}(O_{y_1})$.
Taking into account that $E$
is countable, we can choose an open 
neighborhood $U_1$
of the point $u_1$
such that  $\overline{U_1}\subset\pi^\sharp_{A,K}(O_{y_1}),
\ Fr(U_1)\cap E=\emptyset$
and    $diam(\overline{U_1})<\min\{cs,1\}$.
As $U_1$
we can take an open 
ball $B(u_1,\rho)$
centered at $u_1$
with a sufficiently small radius $\rho\not=d(u_1,u_k),\ k\in N$. 
If $y\in Fr(U_1)$ then $d(u_1,y)=\rho$. Therefore $Fr(U_1) \cap E=\emptyset$.

Assume that the sets $U_k$
are already constructed for all $k < n$ 
such that their closures are disjoint,            
$\overline{U_k}\subset\pi^\sharp_{A,K}(O_{y_k}),\ Fr(U_k)\cap E=\emptyset,\
diam(\overline{U_k})<\min\{cs,1/k\}$    
and   $u_m\in\cup_{k<n}U_k$ for $m<n$.

Step $ n $.  If $ u_n\in \cup_{k <n} U_k $, then put $ U_n = \emptyset $.  Otherwise $ u_n\not \in \cup_{k <n} \overline{U_k} $ and 
 we construct a neighborhood $ U_n $ of the point $ u_n$ as in step 1 such that
$\overline{U_n}\subset \pi_{A,K}^\sharp (O_{y_n})\setminus \cup_{k<n}\overline{U_k}$, $Fr(U_n)\cap E=\emptyset$, $diam \overline{U_n}<\min\{cs,1/n\}$. $\Box$

Set $G=Y_{A,K}\setminus\cup_{n\in N}U_n$
and $W=\pi_{A,K}^{-1}(\cup_{n\in N}U_n)$.
The map  $\pi_{A,K}$ is one-to-one on $Y_A\setminus W$
and $\pi_{A,K}(Y_A\setminus W)=G$.
Define a real function on $G$
 by formula
$$
r=g\circ(\pi_{A,K}^{-1}|_G). \eqno(3.8)
$$
Set
$$
a_n= \inf(g(O_{y_n})),\ b_n=\sup(g(O_{y_n})).\eqno(3.9)
$$
We have
$$ 
b_n-a_n=diam(g(O_{y_n}))<cs.
$$
Thus $r(Fr(U_n))\subset[a_n,b_n]$ since $\overline{U_n}\subset\pi^\sharp_{A,K}(O_{y_n})$.

Let us mention that from Proposition 3.1 and (3.6), (3.9) we 
get
$$
\lim_n(b_n-a_n)=0. \eqno(3.10)
$$
By the Tietze theorem on the extension of continuous functions defined on a closed subset of a metric space, there is a continuous function $ r_n:\overline{U_n} \to [a_n, b_n] $, which is an extension
of the  function $ r|_{Fr (U_n)} $.  

{\bf Claim 2.} {\it The real function                   
$$
q(x)=\left\{
\begin{array}{ll}
r(x)\ if\ x\in G\\
r_n(x)\ if\ x\in U_n,\ n\in N
\end{array}\right.
$$
is continuous on $Y_{A,K}$.}

{\bf Proof of Claim 2.}  Pick  $x\in Y_{A,K}$.
If $ x \in \cup U_n $, then the continuity of $ q$ at the point $ x $ is obvious.

 Let $ x \in G \setminus \cup \overline{U_n} $. Let $z$
 be the unique point of $Y_A$
for which $\pi_{A,K}(z)=x$.
Fix $\varepsilon>0$.
By the continuity of g, 
there exists a neighborhood $V_z$
of $z$ such that for any $v\in V_z$
we have $|g(z)-g(v)|<\varepsilon.$
It is clear that $x\in\pi_{A,K}^\sharp(V_z)$.        
Choose $k\in N$
such that   $B(x,1/k)\subset\pi^\sharp_{A,K}(V_z)$.
Using  (3.10) 
 we may assume in addition  
$$
b_n-a_n<\varepsilon\ for\ n>k. \eqno(3.11)
$$  
We show that $|q(x)-q(w)|<2\varepsilon$ 
for any $w\in B(x,1/2k)\setminus \cup_{n\leq 2k}\overline{U_n}$ (the last set is an open neighborhood of the point $x$).

If $w\in G$,
then there is a unique point $v\in V_z$
such that  $\pi_{A,K}(v)=w$
and $q(w)=g(v)$.
Consequently, 
$$
|q(x)-q(w)|=|g(z)-g(v)|<\varepsilon.
$$
If $w\not\in G$, 
then  $w\in U_n$
for some $n>2k$. Since $diam(U_n)<1/n\leq 1/2k$ and $w\in B(x,1/2k)$,
we get $U_n\subset B(x,1/k)\subset\pi^\sharp_{A,K}(V_z)$.
By virtue of (3.11) we have
$$
diam(g(O_{y_n}))=b_n-a_n<\varepsilon.
$$
Let $t\in\pi_{A,K}^{-1}(w)$.
Then $t\in O_{y_n}$,
since  by construction $U_n\subset\pi^\sharp_{A,K}(O_{y_n})$. Therefore,
$g(t)\in[a_n,b_n]$. By definition $q(w)=r_n(w)\in[a_n,b_n]$.
Thus, $|g(t)-q(w)|<\varepsilon.$
Since $t\in V_z$
 we have
$$
|q(x)-q(w)|=|g(z)-q(w)|\leq|g(z)-g(t)|+|g(t)-q(w)|<2\varepsilon.
$$
The continuity of $q$ at the point $x$ is proved.
 
It remains to consider the case when $ x \in Fr (U_n) $ for some $ n $.  In this case, we represent $ Y_{A,K} $ as a union of two closed sets:
 $$
Y_{A,K}=\overline{U_n}\cup (Y_{A,K}\setminus U_n).
$$
The restriction of $ q$ to $\overline{U_n} $ is continuous at the point $ x $ by the definition of $ q $.  The continuity of the restriction of $ q$ to $ (Y_{A,K} \setminus U_n) $ at the point $ x $ is proved by a verbatim repetition of the above arguments with $ Y_{A,K} $ replaced by $ (Y_{A,K} \setminus U_n) $, $ Y_A $  by $ Y_A \setminus \pi_{A,K}^{- 1}(U_n) $, and $ U_n $  by the empty set.  The continuity of $ q$ at the point $ x $ on $\overline{U_n} $ and $ (Y_{A,K} \setminus U_n) $ implies continuity at the point $ x $ on $ Y_{A,K} $. $\Box$

Set $p=q\circ \pi_{A,K}\circ f_A$.
We have $p\in C(X)$.

{\bf Claim 3.} {\it We have

$(a)\ ||p-h||_\infty\leq c||\Omega_f(h)|_{A\setminus K}||_\infty;$

$(b)\ supp\ \Omega_f(p)\subset K\subset A.$}

{\bf Proof of Claim 3.} (a) Pick $x\in X$. Assume that 
 and $\pi_{A,K} (f_A(x))\in G$. From the definition of $q$ we get $p(x)=h(x)$.

Now let $\pi_{A,K} (f_A(x))\not\in G$. In this case $y=\pi_{A,K} (f_A(x))\in U_n$ for some $n$.   Since $U_n\subset \pi_{A,K}^\sharp(O_{y_n}$)) we get $f_A(x)\in O_{y_n}$
and $h(x)=g(f_A(x))\in[a_n,b_n]$
by virtue 
of (3.9). 
 
On the other hand, when $y\in U_n$
$$
p(x)=q(y)=r_n(y)\in[a_n,b_n],
$$
therefore $|p(x)-h(x)|\leq b_n-a_n\leq cs=c||\Omega_f(h)|_{A\setminus K}||_\infty$. Thus $||p-h||_\infty\leq c||\Omega_f(h)|_{A\setminus K}||_\infty$.

(b) If $y\not\in K\subset Y$, then $\varkappa^{-1}_{A,K}(y)=\{t\}$ is a singleton. Hence for any $x\in f^{-1}(y)$ $p(x)=q(t)$. Therefore $osc_{f^{-1}(y)}(p)=0$ for $y\not\in K$. Thus  $supp \ \Omega_f(p)\subset K$. $\Box$

Since $c$ is arbitrary number bigger than $1$, (a) and (b) 
complete the proof of proposition.
$\Box$

In order to finish the proof of Theorem 1.4 it is 
enough to prove the following

{\bf Proposition 3.5.} {\it  The map of the fiberwise oscillation 
$\Omega_f:C(X)\to c_0(Y)$ satisfies the conditions (a) and (b)
of (i) Theorem 1.1  for any fully closed map  $f:X\to Y$ of compact $X$   
onto a metric compact $Y$ with metrizable fibers $f^{-1}(y),\ y\in Y$.}

{\bf Proof.} Given  a subset $A$ of $Y$
we define $Z_A$ as in (3.1).  According Lemma3.3  $Z_A$
can be identify 
with $C(Y_A)$.
 If $A$ is finite or countable  we get that 
$Y_A$ 
is metrizable by Proposition 2.3. Hence $C(Y_A)$
 is separable. 
Evidently condition (a) holds. Pick $h\in C(X)$
 and set
 $$
K_m=\{y\in Y:\Omega_f(h)(y)\geq 1/m\}.
$$
Applying Proposition 3.1 we get that $K_m,\ m\in N$
are finite. Clearly that $\{K_m:m\in N\}$ is an increasing sequence.
From (3.5) we get
$dist(h,Z_{K_m})\leq 1/m.$
So $\lim_m\ dist(h,Z_{K_m})=0.$ 
From this it follows that
$$
h\in \overline{span\cup_{m\in N}Z_{K_m}}^{||\cdot||_\infty}.
$$
So (b’) holds. 
$\Box$

\begin{center}
{\bf References}
\end{center}

[1] G.Alexandrov,  Spaces of continuous functions isomorphic to 
locally uniformly rotund Banach space, C. R. Acad. Bulg. Sci., 41,    
no 8 (1988), 9-12.

[2]  P.S.Alexandrov, B.A.Pasynkov. Introduction into Dimension Theory. Moscow, "Nauka", 1973 [in Russian]. 

[3] R.Deville, G.Godefroy, V.,Zizler,  Smoothness and renorming 
in Banach spaces, Pitman Monographs and Surveys in Pure and Appl. 
Math. 64, Longman Scientific and Technical, Longman House, Burnt 
Mill, Harlow, 1993.

[4] R.Engelking,  On functions defined in Cartesian products,     
Fundamenta Math., 59 (1966), 221–231.

[5] V.V. Fedorchuk, Fully closed mappings and their applications, J. 
Math. Sci., 136(5) (2006), 4201– 4
291.

[6] S.P.Gul'ko, A.V.Ivanov.  On fully closed mappings of  Fedorchuk compacta.  Vestn. Tom. gos. un-ta. Matematika i mekhanika. 2017. No 50. P. 5-8.

[7] R.Haydon, Trees in renorming theory, Proc. London Math. Soc., 
78 (1999), 541–584.  

[8] R.Haydon, Locally uniformly rotund norms in Banach spaces 
and their duals. J. Funct. Analysis, 254 (2008), 2023–2039. 

[9] R.Haydon, J. Jayne, I. Namioka, C.A. Rogers,  Continuous 
functions on totally ordered spaces that are compact in their order 
topologies. J. Funct. Analysis, 178 (2000), 23–63.  

[10] R. Haydon, A.Molt\'o, J. Orihuela, Spaces of functions with 
countably many discontinuities, Israel J. Math., 158 (2007), 19–39.

[11] A.V. Ivanov, On Fedorchuk compacta, in: Mappings and Functors, 
Izd. Mosk. Univ.1984,pp. 31-40. (in Russian)

[12] J. F. Mart\'inez, A.Molt\'o, J. Orihuela, S.Troyanski, On locally 
uniformly rotund renormings in C(K) spaces, Canadian J. Math., 62 
(2010), 595 - 613. 

[13]  A.Molt\'o, J.Orihuela and S. Troyanski, Locally uniformly 
rotund renorming and fragmentability, Proc. London 
Math. Soc.,75 (1997), 619-640.

[14] A.Molt\'o, J. Orihuela, S.Troyanski, M. Valdivia, Non-linear 
transfer technique, Lect. Notes Math.,1951, Springer, Berlin, 2009.

\end{document}